\documentclass[a4paper,12 pt]{article}
\usepackage{amsmath,amsthm,amssymb,latexsym}
\title{On representations of variants of semigroups}
\author{Ganna Kudryavtseva and Victor Maltcev}
\date{}
\begin{document}

\maketitle

\begin{abstract}
We construct a family of representations of an
arbitrary variant $S_a$ of a semigroup $S$, induced by a given
representation of $S$, and investigate properties of such
representations and their kernels.
\end{abstract}

\newtheorem{theorem}{Theorem}
\newtheorem{lemma}{Lemma}
\newtheorem{proposition}{Proposition}
\newtheorem{corollary}{Corollary}
\newtheorem{definition}{Definition}
\newtheorem{example}{Example}
\newtheorem{remark}{Remark}

2000 Mathematics Subject Classification 20M10, 20M30.

\section{Introduction} \label{sec:intro}
Let $S$ be a semigroup and $a\in S$. Set $s\ast t=sat$,
$s,t\in S$. The semigroup $(S,\ast)$  is called a
\emph{variant} of $S$ and will be denoted by $S_a$. One of the motivations for the study of
variants is their importance to
understanding the structure of an original semigroup. In
particular, the notion of a regularity-preserving element, which
generalizes the notion of an invertible element, is defined using
variants, see \cite{Hic1}, \cite{Hic2}, \cite{Kh_Laws}. Other
situations where variants naturally appear and work are discussed
in \cite{Kh_Laws}.

Let $\cdot: M\times S\to M$ or just $(S,M)$ be a representation of
a semigroup $S$ by transformations of a set $M$.
It
defines a homomorphism from $S$ to the full transformation
semigroup $\mathcal{T}(M)$; we denote the
image of $(m,s)$ under the function $\cdot$ by $m\cdot s$.
For  $a\in S$ and a decomposition $a=\beta\alpha$,
$\alpha,\beta\in S^{1}$, we introduce a map $\ast: M\times S_a\to
M$ or just $(S_a,M;\alpha,\beta)$, defined by
\begin{equation*}
m\ast s=m\cdot(\alpha s\beta)~\mbox{for all}~s\in S,~m\in M.
\end{equation*}
 Since \begin{multline*}
m\ast (s\ast t)=m\ast (sat)=m\cdot (\alpha s\beta\alpha
t\beta)=\\(m\cdot (\alpha s\beta))\cdot (\alpha t\beta)=(m\ast
s)\ast t$ for all $s,t\in S,~m\in M,
\end{multline*} it follows that
 $(S_a,M;\alpha,\beta)$ is a representation of $S_a$ on $M$.

The kernel $\rho$ of a representation of $S$
naturally leads to the consideration of the family of congruences
$\{\rho_{b,c}:b,c\in S^{1},cb=a\}$ on $S_a$, which are the kernels of the corresponding
representations of $S_a$.
In the case when $\iota$ is the identity relation on $S$, $\iota_{a,1}$ and
$\iota_{1,a}$ coincide with the congruences $\mathrm{l}$ and
$\mathrm{r}$ on $S_a$ respectively, which were first introduced by
Symons (see \cite{Symons}) for the generalized
transformation semigroups and many properties of which resemble
the corresponding properties of the Green's relations
$\mathcal{L}$ and $\mathcal{R}$.

The present paper is devoted to the study of the connections between
a given representation $(S,M)$ of
$S$ on $M$ and the representations $(S_a,M;\alpha,\beta)$; and
also to the study of the properties of congruences $\{\rho_{b,c}:b,c\in S^{1},cb=a\}$.
In Section \ref{sec:2} we prove that for a regular semigroup $S$, which does not contain the bicyclic semigroup, a faithful
representation of $S_a$ coincides with certain $(S_a,M;\alpha,\beta)$ if and only if
$\alpha$ and $\beta$ are left and right cancellable in $S$ respectively (Theorems \ref{th:faith_rep},
\ref{th:reg} and Proposition \ref{aaa}).
We also show that for the bicyclic semigroup these conditions are
not equivalent. The next result is Theorem
\ref{th:easy}, claiming that for a regular semigroup $S$ each representation
$\ast:~M\times S_a\rightarrow M$ of any variant of $S$ is induced by some
$\cdot: M\times S\rightarrow M$ such that either $m\ast s=m\cdot
(as)$ for all $m\in M$ or $m\ast s=m\cdot (sa)$ for all $m\in M$,
if and only if it is either a left or right group, giving a new characterization of the class
of left (right) groups. In Section
\ref{sec:3} we  investigate
the  properties of the family of congruences
$\{\rho_{b,c}:b,c\in S^{1},cb=a\}$.
 In particular, in Theorem \ref{th:beautiful}
we prove that in the case when $c$, $b$ are regular elements of
$S$ the map sending $\rho$ to $\rho_{b,c}$ is a homomorphism from
the lattice of congruences of $S$ to the lattice of congruences of
$S_{cb}$. Without the requirement of regularity of $b,c$ this map
always preserves the meet, while may not preserve the join, the
corresponding example is provided. Finally, the aim of Section
\ref{sec:4} is, given $a\in S^1$ and a congruence $\rho$ on $S$,
we define a congruence $\rho_a$ which naturally generalizes the
congruence $\mathrm{d}$ on $S_a$, introduced in \cite{Symons} and
studied afterwards in \cite{Hic2} and  \cite{Magill}.

\section{Connections between representations} \label{sec:2}

Let $S$ be a semigroup. An element $u\in S^{1}$ will be called
\emph{left} {\rm [}\emph{right}{\rm ]} \emph{cancellable} if
$us=ut$ {\rm[}$su=tu${\rm]} implies $s=t$ for all $s,t\in S$.

Fix $a\in S$ and a decomposition $a=\beta\alpha$, $\alpha,\beta\in
S^1$.

\begin{lemma} \label{lm:good}
The implication $\alpha s\beta=\alpha
t\beta\Rightarrow~s=t~\mbox{for all}~s,t\in S$ holds if and only
if $\alpha$ is left cancellable and $\beta$ is right cancellable
simultaneously.
\end{lemma}

\begin{proof}
The proof is straightforward.
\end{proof}

A representation $\cdot:M\times S\to M$ is said to be
\emph{faithful} if the corresponding homomorphism from $S$ to
$\mathcal{T}(M)$ is injective. That is, $\cdot:M\times S\to M$ is
faithful if and only if
\begin{equation*}
(\forall m\in M)~[m\cdot s=m\cdot t]\Rightarrow~[s=t]~\mbox{for
all}~s,t\in S.
\end{equation*}

\begin{lemma} \label{lm:useful}
Suppose that $(S_a,M;\alpha,\beta)$ is faithful. Then so is
$(S,M)$.
\end{lemma}

\begin{proof}
Indeed, by the definition we have
\begin{equation*}
(\forall~m\in M)~[m\cdot(\alpha s\beta)=m\cdot(\alpha
t\beta)]\Rightarrow~[s=t]~\mbox{for all}~s,t\in S.
\end{equation*}
Take $s,t\in S$. If now $m\cdot s=m\cdot t$ for all $m\in M$, then
\begin{equation*}
(m\cdot\alpha)\cdot (s\beta)=(m\cdot\alpha)\cdot
(t\beta)~\mbox{for all}~m\in M,
\end{equation*}
and so $s=t$.
\end{proof}

\begin{theorem} \label{th:faith_rep}
Let $(S,M)$ be a representation of a semigroup $S$ on a set $M$
and $a\in S$. Take $\alpha,\beta\in S^{1}$ such that
$\beta\alpha=a$. Then the following conditions are equivalent:
\begin{enumerate}
\item
$(S_a,M;\alpha,\beta)$ is faithful;
\item
$(S,M)$ is faithful and $\alpha,\beta$ are left and right
cancellable respectively.
\end{enumerate}
\end{theorem}

\begin{proof}
That 1) implies 2) is due to Lemmas~\ref{lm:good} and
\ref{lm:useful}. The opposite implication is a consequence of
Lemma~\ref{lm:good}.
\end{proof}

Recall that an element $u$ of a semigroup $S$ is said to be a
\emph{mididentity} in $S$ if $sut=st$ for all $s,t\in S$. For the
case when $S$ is regular, we obtain the following theorem.

\begin{theorem} \label{th:reg}
Let $S$ be a regular semigroup and $a=\beta\alpha\in S$,
$\alpha,\beta\in S^{1}$. Suppose that there exist $\alpha^{\ast}$
and $\beta^{\ast}$ of $S^1$, inverses of $\alpha$ and $\beta$ in
$S^1$ respectively, such that
$\beta^{\ast}\beta\alpha\alpha^{\ast}$ is a mididentity in $S$.
Let also $\ast:M\times S_a\to M$ be a faithful representation of
$S_a$ on $M$. Then the following conditions are equivalent:
\begin{enumerate}
\item
there is a representation $(S,M)$ such that $(S_a,M)$ coincides
with $(S_a,M;\alpha,\beta)$;
\item
$\alpha$ and $\beta$ are left and right cancellable respectively.
\end{enumerate}
\end{theorem}

\begin{proof}
Due to Theorem \ref{th:faith_rep} we have that 1) implies 2).

Now let $\alpha$ and $\beta$ be left and right cancellable
respectively. Let us prove that a function $\cdot:M\times
S\rightarrow M$, given by
\begin{equation*}
m\cdot s=m\ast (\alpha^{\ast}s\beta^{\ast})~\mbox{for all}~s\in
S,~m\in M,
\end{equation*}
defines a representation $(S,M)$ of $S$ such that $(S_a,M)$
coincides with $(S_a,M;\alpha,\beta)$.

Since $\beta^{\ast}\beta\alpha\alpha^{\ast}$ is a mididentity in
$S$, we have
\begin{multline*}
(m\cdot s)\cdot t=(m\ast(\alpha^{\ast}s\beta^{\ast}))\ast
(\alpha^{\ast}t\beta^{\ast})=m\ast
(\alpha^{\ast}s\beta^{\ast}\cdot \beta\alpha\cdot
\alpha^{\ast}t\beta^{\ast})=\\m\ast
(\alpha^{\ast}st\beta^{\ast})=m\cdot(st)~\mbox{for all}~s,t\in
S,~m\in M.
\end{multline*}
Thus $(S,M)$ is indeed a representation.

Note that the equality $\alpha\alpha^{\ast}\alpha s=\alpha s$
implies $\alpha^{\ast}\alpha s=s$ for all $s\in S$. Analogously,
$s\beta\beta^{\ast}=s$ for all $s\in S$. Then we have
\begin{equation*}
m\cdot (\alpha s\beta)=m\ast (\alpha^{\ast}\alpha
s\beta\beta^{\ast})=m\ast s~\mbox{for all}~s\in S,~m\in M.
\end{equation*}
The latter means that $(S_a,M)$ coincides with
$(S_a,M;\alpha,\beta)$. That is, the function $\cdot$ defines the
required representation $(S,M)$.
\end{proof}

The following example shows that the condition of Theorem
\ref{th:reg} that there are $\alpha^{\ast}$ and $\beta^{\ast}$ in
$S^{1}$, such that $\beta^{\ast}\beta\alpha\alpha^{\ast}$ is a
mididentity in $S$, is essential.

\begin{example}
Let $\mathcal{B}=\langle a,b~|~ba=1\rangle$ be the bicyclic
semigroup. Consider the Cayley representation of $\mathcal{B}$:
\begin{equation*}
m\ast s=ms~\mbox{for all}~m,s\in \mathcal{B}.
\end{equation*}
Then there is no faithful representation $\circ: \mathcal{B}\times
\mathcal{B}\rightarrow \mathcal{B}$ such that
\begin{equation*}
m\circ (asb)=m\ast s=ms~\mbox{for all}~m,s\in \mathcal{B}.
\end{equation*}
\end{example}

Since $\mathcal{B}$ is an inverse semigroup with the identity element, $1$,
its Cayley representation is faithful. That
$ba=1$ implies $\mathcal{B}_{ba}=\mathcal{B}$. Note that
$b^{-1}baa^{-1}=a\cdot ba\cdot b=ab$ is not a mididentity in
$\mathcal{B}$, as $1\cdot ab\cdot 1=ab\ne 1$.

Assume that there is a faithful representation $\circ:
\mathcal{B}\times \mathcal{B}\rightarrow \mathcal{B}$ such that
\begin{equation*}
m\circ (asb)=m\ast s=ms~\mbox{for all}~m,s\in \mathcal{B}.
\end{equation*}
Then $m\circ (ab)=m\ast 1=m$ for all $m\in \mathcal{B}$. So we
have
\begin{multline*}
m\circ a=(m\circ a)\circ(ab)=m\circ (a^2b)=m\ast a=ma~\mbox{and}\\
m\circ b=(m\circ(ab))\circ b=m\circ(ab^2)=m\ast b=mb~\mbox{for
all}~m\in \mathcal{B}.
\end{multline*}
But this leads to a contradiction:
\begin{equation*}
1=1\circ(ab)=(1\circ a)\circ b=a\circ b=ab.
\end{equation*}

Let us recall (see Corollary 1.32 of \cite{Clifford}) that every
homomorphic image of $\mathcal{B}$ is either a cyclic group or
isomorphic to $\mathcal{B}$. Using this fact, we are going to
prove that in the case when $S$ is a regular semigroup, which does
not contain a subsemigroup, isomorphic to $\mathcal{B}$, one
obtains that the conditions 1) and 2) of Theorem \ref{th:reg} are
equivalent.

\begin{proposition}\label{aaa}
Let $S$ be a regular semigroup, which does not contain a
subsemigroup isomorphic to $\mathcal{B}$. Let $a=\beta\alpha\in
S$, $\alpha,\beta\in S^{1}$. Let also $\ast:M\times S_a\to M$ be a
faithful representation of $S_a$ on $M$. Then the conditions 1)
and 2) of Theorem \ref{th:reg} are equivalent.
\end{proposition}

\begin{proof}
We have only to prove that 2) implies 1).

Let now $\alpha$ and $\beta$ be left and right cancellable
respectively. Take $\alpha^{\ast}$ and $\beta^{\ast}$, inverses of
$\alpha$ and $\beta$ in $S^1$ respectively. Then the equality
$\alpha\alpha^{\ast}\alpha s=\alpha s$ implies
$\alpha^{\ast}\alpha s=s$ for all $s\in S$. Also
$s\beta\beta^{\ast}\beta=s\beta$ implies $s\beta\beta^{\ast}=s$
for all $s\in S$. It follows that $\alpha^{\ast}\alpha$ and
$\beta\beta^{\ast}$ are left and right identity elements of $S^1$
respectively. Consider three possible cases.

\emph{Case 1}. Let $\alpha,\beta\in S$. Then $\alpha^{\ast}\alpha$
and $\beta\beta^{\ast}$ are both in $S$. It follows that $S$ has
the identity $1=\alpha^{\ast}\alpha=\beta\beta^{\ast}$. Assume
that $\alpha\alpha^{\ast}\ne 1$. Then
$\langle\alpha,\alpha^{\ast}\rangle$ is a homomorphic image of
$\mathcal{B}=\langle a,b\mid ba=1\rangle$ under a map, which sends
$a$ to $\alpha$ and $b$ to $\alpha^{\ast}$. Since
$\langle\alpha,\alpha^{\ast}\rangle$ is not isomorphic to a cyclic
group, then $\langle\alpha,\alpha^{\ast}\rangle$ is a
subsemigroup, isomorphic to $\mathcal{B}$. We get a contradiction.
Thus $\alpha\alpha^{\ast}=1$. Analogously, one shows that
$\beta^{\ast}\beta=1$. Then $\beta^{\ast}\beta\alpha\alpha^{\ast}$
is a mididentity in $S$. Now due to Theorem \ref{th:reg}, we have
that 1) holds.

\emph{Case 2}. Let $\beta\notin S$ and $\alpha=a\in S$. Arguments
are similar to those of Case 1. Then
$\alpha\alpha^{\ast}=\alpha^{\ast}\alpha$ is again a mididentity
in $S$ and to prove that 1) holds, it remains to use Theorem
\ref{th:reg}.

\emph{Case 3}. Let $\alpha\notin S$ and $\beta=a\in S$.
Analogously to Case 2, $\beta^{\ast}\beta=\beta\beta^{\ast}$ is a
mididentity in $S$. Thus, according to Theorem \ref{th:reg}, 1)
holds.
\end{proof}

Recall (see \cite{Clifford}) that a \emph{right} [\emph{left}]
\emph{group} is a semigroup which is right [left] simple and left
[right] cancellative. We will need the following known fact.

\begin{theorem}[Theorem 1.27 of \cite{Clifford}] \label{clifford}
The following conditions for a semigroup $S$ are equivalent:
\begin{enumerate}
\item
$S$ is a left {\rm[}right{\rm]} group;
\item
$S$ is right {\rm[}left{\rm]} simple and contains an idempotent;
\item
$S$ is a direct product $G\times U$ of a group $G$ and a left
{\rm[}right{\rm]} zero semigroup $U$.
\end{enumerate}
\end{theorem}

Among all the decompositions $a=\beta\alpha$, $\alpha,\beta\in
S^{1}$, two are rather special, namely $a=a\cdot 1$ and $a=1\cdot
a$. The main result of the following theorem is the
characterization of semigroups $S$ such that for all $a\in S$ and
for each representation $\ast:~M\times S_a\rightarrow M$ there
exists a representation $\cdot: M\times S\rightarrow M$ such that
either $m\ast s=m\cdot (as)$ for all $m\in M$ or $m\ast s=m\cdot
(sa)$ for all $m\in M$.

\begin{theorem} \label{th:easy}
Let $S$ be a regular semigroup with the set of idempotents $E$.
Then the following conditions are equivalent:
\begin{enumerate}
\item
for all $a\in S$ and for each representation $\ast:~M\times
S_a\rightarrow M$ there exists a representation $\cdot: M\times
S\rightarrow M$ such that either $m\ast s=m\cdot (as)$ for all
$m\in M$ or $m\ast s=m\cdot (sa)$ for all $m\in M$;
\item
for all $e\in E$ there are a faithful representation
$\ast:~M\times S_e\rightarrow M$ and a representation $\cdot:
M\times S\rightarrow M$ such that either $m\ast s=m\cdot (es)$ for
all $m\in M$ or $m\ast s=m\cdot (se)$ for all $m\in M$;
\item
every element of $E$ is either a left or a right identity  in $S$;
\item
$S$ is either a left or a right group.
\end{enumerate}
\end{theorem}

\begin{proof}
Obviously, 1) implies 2).

Assume that 2) holds. Take an idempotent $e\in E$. The first case
is that there are a set $M$, a faithful representation
$\ast:~M\times S_e\rightarrow M$ and a representation $\cdot:
M\times S\rightarrow M$ such that $m\ast s=m\cdot (es)$ for all
$m\in M$. Then $m\ast s=m\ast (es)$ for all $m\in M$ and $s\in S$
and so $s=es$ for all $s\in S$. In the second case we obtain
$s=se$ for all $s\in S$. Thus, 2) implies 3).

Let us prove now that 3) implies 4). Note that $E$ is nonempty as
$S$ is regular. If $S$ contains both left and right identities
then it contains the identity $1$, whence each idempotent
coincides with $1$. Hence, $S$ is a regular semigroup with the
unique idempotent $1$. It follows now from Lemma II.2.10 from
\cite{Petrich} that $S$ is a group and so also a left group. If
every idempotent $e\in E$ is a left {\rm[}right{\rm]} unit in $S$,
then $sS^{1}=sS=ss'S=S$ {\rm[}$S^{1}s=Ss=Ss's=S${\rm]} for all
$s\in S$ and for each $s'$, inverse of $s$, and so $S$ is right
{\rm[}left{\rm]} simple. But due to Theorem \ref{clifford} we
obtain now that $S$ is a right [left] group.

Finally, assume that $S$ is, for instance, a left group. Then by
Theorem \ref{clifford} there are a group $G$ and a left zero
semigroup $U$ such that $S=G\times U$. Take $(g,u)\in S$ and a
representation $\ast: M\times S_{(g,u)}\rightarrow M$. Put $m\cdot
(h,v)=m\ast (hg^{-1},v)$ for all $(h,v)\in S$. Then $m\cdot
((h,v)(g,u))=m\ast (h,v)$ for all $(h,v)\in S$. We are left to
prove that $\cdot$ is a representation of $S$. Indeed, we have
\begin{multline*}
m\cdot (h_1h_2,v_1v_2)=m\ast (h_1h_2g^{-1},v_1v_2)=m\ast
((h_1g^{-1},v_1)\ast(h_2g^{-1},v_2))=\\(m\cdot(h_1,v_1))\cdot
(h_2,v_2)~\mbox{for all}~(h_1,v_1),(h_2,v_2)\in S,~m\in M.
\end{multline*}
This completes the proof.
\end{proof}

Let $S$ be a semigroup and $\cdot:M\times S\to M$ a representation
of $S$. Denote by $M\cdot\alpha$ the set $\{m\cdot\alpha:~m\in M
\}$ and by $m\cdot S$ the set $\{m\cdot s:~s\in S\}$. In the case
when $S\ne S^{1}$ set $M\cdot 1=M$. A representation $(S,M)$ is
said to be \emph{cyclic} (see \cite{Clifford}) if there is a
\emph{generating} element $m\in M$ for $(S,M)$, i.e., $m\cdot
S=M$.

\begin{proposition}
Let $S$ be an arbitrary semigroup and $a\in S$, $a=\beta\alpha$,
$\alpha,\beta\in S^{1}$. Let $(S,M)$ be a representation of $S$ on
$M$. The following conditions are equivalent:
\begin{enumerate}
\item
$(S_a,M;\alpha,\beta)$ is cyclic;
\item
$(S,M)$ is cyclic, $M\cdot\beta=M$ and $M\cdot\alpha$ contains a
generating element for $(S,M)$.
\end{enumerate}
\end{proposition}

\begin{proof}
Suppose that 1) holds. It follows immediately from the definition
of cyclic representations that $(S,M)$ is cyclic. If $m\in M$ is a
generating element for $(S_a,M;\alpha,\beta)$ then $m\cdot (\alpha
S\beta)=M$ which implies $M\cdot\beta=M$. Since
$(m\cdot\alpha)\cdot S\supseteq m\cdot(\alpha S\beta)=M$, we have
that $M\cdot\alpha$ contains a generating element for $(S,M)$.

Suppose now that 2) holds. Let $m\in M\cdot\alpha$ be a generating
element for $(S,M)$ such that $m=m_0\alpha$ for some $m_0\in M$.
Let us prove that $m_0$ is a generating element for
$(S_a,M;\alpha,\beta)$. Indeed, take $m_1\in M$. Then there exists
$\overline{m}\in M$ such that $m_1=\overline{m}\cdot\beta$. Also
there is $s\in S$ such that $\overline{m}=(m_0\alpha)\cdot s$.
Then we obtain $m_0\ast s=m_0\cdot (\alpha
s\beta)=((m_0\alpha)\cdot s)\cdot\beta=m_1$. Thus, $m_0\ast S=M$.
\end{proof}

\section{Congruences $\rho_{b,c}$} \label{sec:3}

Let $\rho$ be a congruence on a semigroup $S$ and $a\in S$,
$a=cb$, $b,c\in S^1$. Define a relation $\rho_{b,c}$ on $S_{a}$ as
follows:
\begin{equation*}
(s,t)\in\rho_{b,c}~\mbox{if and only
if}~(bsc,btc)\in\rho~\mbox{for all}~s,t\in S.
\end{equation*}
It is straightforward that $\rho_{b,c}$ is a congruence on $S_a$.
If now we have a representation $\cdot:M\times S\rightarrow M$ of
$S$ on $M$, then the congruence $\nu$, related to it, is given by
\begin{equation*}
(s,t)\in\nu~\mbox{if and only if}~m\cdot s=m\cdot t~\mbox{for
all}~m\in M.
\end{equation*}
Then one can easily see that the congruence on $S_a$, related to
the representation $(S_a,M;b,c)$, coincides with $\nu_{b,c}$.

Denote by $\mathrm{Cong}(S)$ the set of all congruences on a
semigroup $S$. Set $\rho_{1,1}=\rho$ for all
$\rho\in\mathrm{Cong}(S)$. In the case when $S\ne S^1$, set also
$S_{1}=S$.

\begin{proposition} \label{pr:fact}
Let $S$ be a semigroup, $b,c\in S^{1}$, $\rho\in\mathrm{Cong}(S)$.
Then
\begin{equation*}
S_{cb}/\rho_{b,c}\cong bSc/\rho\cap(bSc\times bSc).
\end{equation*}
\end{proposition}

\begin{proof}
Define a map $\varphi: S_{cb}\rightarrow bSc/\rho\cap(bSc\times
bSc)$ as follows:
\begin{equation*}
\varphi(x)=(bxc)\rho~\mbox{for all}~x\in S.
\end{equation*}
One can easily show that $\varphi$ is an onto homomorphism and
$\varphi\circ\varphi^{-1}=\rho_{b,c}$. These two facts complete
the proof.
\end{proof}

Let $S\ne S^{1}$ and $\rho\in\mathrm{Cong}(S)$. Then we identify
$(1)\rho$ with the identity element of $(S/\rho)^1$.

\begin{theorem} \label{th:simple}
Let $S$ be a semigroup, $b,c,b_1,c_1\in S^1$,
$\rho,\sigma\in\mathrm{Cong}(S)$. Then
\begin{enumerate}
\item
if $b\rho$ and $b_1\rho$ are $\mathcal{L}$-related, $c\rho$ and
$c_1\rho$ are $\mathcal{R}$-related in $(S/\rho)^1$, then
$\rho_{b,c}=\rho_{b_1,c_1}$;
\item
$\rho\subseteq\rho_{b,c}$;
\item
$\rho_{b,c}=\rho$ if and only if $b\rho$ and $c\rho$ are left and
right cancellable in $(S/\rho)^1$ respectively;
\item
if $\rho\subseteq\sigma$ then $\rho_{b,c}\subseteq\sigma_{b,c}$.
If $b\sigma$ and $c\sigma$ are left and right cancellable in
$(S/\sigma)^1$ then $\rho\subseteq\sigma$ if and only if
$\rho_{b,c}\subseteq\sigma_{b,c}$.
\end{enumerate}
\end{theorem}

\begin{proof}
Statements 1) and 2) follow immediately from the definition of
$\rho_{b,c}$. Statement 3) follows from Lemma \ref{lm:good} and
the fact that $\rho_{b,c}=\rho$ is equivalent to the implication
\begin{equation*}
(b\rho)x(c\rho)=(b\rho)y(c\rho)\Rightarrow [x=y]~\mbox{for
all}~x,y\in S/\rho.
\end{equation*}
Finally, let us prove 4). If $\rho\subseteq\sigma$ then
$(x,y)\in\rho_{b,c}$ implies $(bxc,byc)\in\rho\subseteq\sigma$
which, in turn, implies $(x,y)\in\sigma_{b,c}$. Thus, if
$\rho\subseteq\sigma$, then $\rho_{b,c}\subseteq\sigma_{b,c}$.

Suppose now that $b\sigma$ and $c\sigma$ are left and right
cancellable in $(S/\sigma)^1$ and
$\rho_{b,c}\subseteq\sigma_{b,c}$. Take $(x,y)\in\rho$. Then
$(bxc,byc)\in\rho$ which implies
$(x,y)\in\rho_{b,c}\subseteq\sigma_{b,c}$ or just
$(b\sigma)(x\sigma)(c\sigma)=(b\sigma)(y\sigma)(c\sigma)$, whence
$x\sigma=y\sigma$. Thus, we obtain that $\rho\subseteq\sigma$.
This completes the proof.
\end{proof}

The converse statement of 1) of Theorem \ref{th:simple} is not
true in general as the following easy example shows.

\begin{example}
Consider a semilattice $E=\{a,b\}$, where $a\leq b$. Let $\rho$ be
the identity relation on $E$. Then $\rho_{a,a}=\rho_{a,b}$ but
$(a,b)\notin\mathcal{R}$.
\end{example}

Set $\rho_{a}^r=\rho_{1,a}$ and $\rho_{a}^l=\rho_{a,1}$ for all
congruences $\rho$ on a semigroup $S$ and $a\in S^1$. The
following proposition shows that the converse statement of 1) of
Theorem \ref{th:simple} is true in the case when $b=b_1=1$
[$c=c_1=1$] and $S/\rho$ is inverse.

\begin{proposition}
Let $S$ be a semigroup, $b,c\in S^1$, $\rho\in\mathrm{Cong}(S)$.
Suppose that $S/\rho$ is inverse. Then $\rho_{b}^r=\rho_{c}^r$
{\rm[}$\rho_{b}^l=\rho_{c}^l${\rm]} if and only if $b\rho$ and
$c\rho$ are $\mathcal{R}$-related
{\rm[}$\mathcal{L}$-related{\rm]} in $(S/\rho)^1$.
\end{proposition}

\begin{proof}
The sufficiency follows from 1) of Theorem \ref{th:simple}. Let
now assume that $\rho_{b}^r=\rho_{c}^r$. Set $b\rho=u$ and
$c\rho=v$. Then we have $xu=yu$ if and only if $xv=yv$ for all
$x,y\in S/\rho$. In particular, $xu=xuu^{-1}u$ gives us
$xv=xuu^{-1}v$ for all $x\in S/\rho$. Hence,
$v=vv^{-1}v=vv^{-1}uu^{-1}v=uu^{-1}v$, whence
$vv^{-1}=uu^{-1}vv^{-1}$. Analogously, $uu^{-1}=vv^{-1}uu^{-1}$.
Thus, $vv^{-1}=uu^{-1}$ which is well-known to be equivalent to
$u\mathcal{R}v$ (see \cite{Higgins}).
\end{proof}

Take $b,c\in S^1$. Set a map $\varphi_{b,c}:
\mathrm{Cong}(S)\rightarrow\mathrm{Cong}(S_{cb})$ as follows:
\begin{equation*}
\varphi_{b,c}(\rho)=\rho_{b,c}~\mbox{for
all}~\rho\in\mathrm{Cong}(S).
\end{equation*}
It is well-known that if one has $\rho,\sigma\in\mathrm{Cong}(S)$
then $\rho\vee\sigma$ coincides with the transitive closure
$(\rho\cup\sigma)^{t}$. Denote by $\mathcal{SL}(S)$ the lower
semilattice $(\mathrm{Cong}(S);\subseteq,\cap)$ of congruences on
$S$ and by $\mathcal{L}(S)$ the lattice
$(\mathrm{Cong}(S);\mbox{$\subseteq,$}\cap,\vee)$ of congruences
on $S$.

\begin{theorem} \label{th:beautiful}
Let $S$ be a semigroup and $b,c\in S^1$. Then
\begin{enumerate}
\item
$\varphi_{b,c}$ is a homomorphism from $\mathcal{SL}(S)$ to
$\mathcal{SL}(S_{cb})$;
\item
$\rho_{b,c}\vee\sigma_{b,c}\subseteq(\rho\vee\sigma)_{b,c}$ for
all $\rho,\sigma\in\mathrm{Cong}(S)$;
\item
if $b$ and $c$ are regular in $S^1$ then $\varphi_{b,c}$ is a
homomorphism from $\mathcal{L}(S)$ to $\mathcal{L}(S_{cb})$;
\item
if $bSc=S$ then $\varphi_{b,c}$ is injective homomorphism from
$\mathcal{L}(S)$ to $\mathcal{L}(S_{cb})$.
\end{enumerate}
\end{theorem}

\begin{proof}
Take $\rho,\sigma\in\mathrm{Cong}(S)$. Then
$(x,y)\in\rho_{b,c}\cap\sigma_{b,c}$ if and only if
$(bxc,byc)\in\rho\cap\sigma$, which is, in turn, equivalent to
$(x,y)\in(\rho\cap\sigma)_{b,c}$ for all $x,y\in S$. Thus,
$\rho_{b,c}\cap\sigma_{b,c}=(\rho\cap\sigma)_{b,c}$, which
completes the proof of 1).

To prove 2), we note that due to
$\rho_{b,c}\subseteq(\rho\vee\sigma)_{b,c}$ and
$\sigma_{b,c}\subseteq(\rho\vee\sigma)_{b,c}$, we have that
$\rho_{b,c}\vee\sigma_{b,c}\subseteq(\rho\vee\sigma)_{b,c}$.

Let now $b$ and $c$ be regular in $S^1$ and $b'$, $c'$ be certain
inverses of $b$ and $c$ in $S^1$ respectively. To prove 3), in
view of what has already been done, we are left to show that
$(\rho\vee\sigma)_{b,c}\subseteq\rho_{b,c}\vee\sigma_{b,c}$.
Suppose that $(x,y)\in(\rho\vee\sigma)_{b,c}$. Then
$(bxc,byc)\in\rho\vee\sigma$ and so there are $p_1,\ldots,p_m\in
S$ such that
\begin{equation} \label{eq:music}
(bxc,p_1)\in\rho\cup\sigma,\ldots,(p_{i},p_{i+1})\in\rho\cup\sigma,\ldots,
(p_{m},byc)\in\rho\cup\sigma.
\end{equation}
Clearly, $\rho\cup\sigma$ is left and right compatible, and
therefore
\begin{multline*}
(bxc,b\cdot b'p_1c'\cdot c)\in\rho\cup\sigma,\ldots,\\(b\cdot
b'p_{i}c'\cdot c,b\cdot b'p_{i+1}c'\cdot
c)\in\rho\cup\sigma,\ldots,(b\cdot b'p_{m}c'\cdot
c,byc)\in\rho\cup\sigma,
\end{multline*}
which yields
\begin{multline*}
(x,b'p_1c')\in\rho_{b,c}\cup\sigma_{b,c},\ldots,(b'p_{i}c',b'p_{i+1}c')\in\rho_{b,c}\cup\sigma_
{b,c},\ldots,\\(b'p_{m}c',y)\in\rho_{b,c}\cup\sigma_{b,c},
\end{multline*}
whence $(x,y)\in\rho_{b,c}\vee\sigma_{b,c}$. Thus, 3) is proved.

Finally, assume that $bSc=S$. To prove that $\varphi_{b,c}$ is a
homomorphism, it is enough to show that
$(\rho\vee\sigma)_{b,c}\subseteq\rho_{b,c}\vee\sigma_{b,c}$. Take
$(x,y)\in(\rho\vee\sigma)_{b,c}$. Then there are
$p_1,\ldots,p_m\in S$ such that \eqref{eq:music} holds. Then due
to $bSc=S$ and the fact that $(bsc,btc)\in\rho\cup\sigma$ if and
only if $(s,t)\in\rho_{b,c}\cup\sigma_{b,c}$ for all $s,t\in S$,
we obtain that $(x,y)\in\rho_{b,c}\vee\sigma_{b,c}$. It remains to
prove that $\varphi_{b,c}$ is injective. Suppose that
$\rho_{b,c}=\sigma_{b,c}$. Take $(s,t)\in \rho$. There are $s_1$
and $t_1$ of $S$ such that $s=bs_1c$ and $t=bt_1c$. Then
$(s_1,t_1)\in\rho_{b,c}$ which implies $(s_1,t_1)\in\sigma_{b,c}$,
which, in turn, is equivalent to $(s,t)\in\sigma$. Thus,
$\rho\subseteq\sigma$. Analogously, $\sigma\subseteq\rho$. So
$\rho=\sigma$ and $\varphi_{b,c}$ is an injective homomorphism.
\end{proof}

We note that the converse inclusion of one in 2) of Theorem
\ref{th:beautiful}, namely
$(\rho\vee\sigma)_{b,c}\subseteq\rho_{b,c}\vee\sigma_{b,c}$ for
$\rho,\sigma\in\mathrm{Cong}(S)$, is not true in general as the
following example illustrates.

\begin{example}
Consider the free semigroup $\{a,b\}^{+}$ over the alphabet
$\{a,b\}$. Let $I$ be the ideal consisting of all words from $\{a,b\}^{+}$ of
length not less than $3$. Set $S=\{a,b\}^{+}/I$ to be the Rees quotient semigroup.
\end{example}

Set $\rho=(ba,ab)\cup(ab,ba)\cup\iota$ and
$\sigma=(ab,bb)\cup(bb,ab)\cup\iota$, where $\iota$ denotes the
identity relation on $S$. It follows immediately from the
construction of $S$ that $\rho,\sigma\in\mathrm{Cong}(S)$. Now we
observe that $(a,b)\in(\rho\vee\sigma)_{b,1}$. Indeed, we have
$(ba,ab)\in\rho$ and $(ab,bb)\in\sigma$, whence
$(ba,bb)\in\rho\vee\sigma$ which is equivalent to
$(a,b)\in(\rho\vee\sigma)_{b,1}$. But
$(a,b)\notin\rho_{b,1}\vee\sigma_{b,1}$. Indeed, otherwise there
would exist $t_1,\ldots,t_n\in S$ such that
\begin{equation*}
(ba,bt_1)\in\rho\cup\sigma,\ldots,(bt_{i},bt_{i+1})\in\rho\cup\sigma,\ldots,
(bt_{n},bb)\in\rho\cup\sigma.
\end{equation*}
It follows from the construction of $\rho$ and $\sigma$ that
$(ba,bt_1)\in\rho\cup\sigma$ implies $t_1=a$. Now inductive
arguments show that $t_i=a$ for all possible $i$. In particular,
$t_n=a$. Then $(ba,bb)\in\rho\cup\sigma$, and we get a
contradiction. Thus,
$(\rho\vee\sigma)_{b,1}\nsubseteq\rho_{b,1}\vee\sigma_{b,1}$.

\section{Congruences $\rho_{a}$} \label{sec:4}

Let now $\rho$ be a congruence on a semigroup $S$ and $a\in S^1$.
Then we can construct a congruence $\rho_a\in\mathrm{Cong}(S_a)$
as follows:
\begin{equation*}
(s,t)\in\rho_a~\mbox{if and only if}~(asa,ata)\in\rho~\mbox{for
all}~s,t\in S.
\end{equation*}
Thus, in terms of \cite{Hic2}, $(s,t)\in\rho_a$ if and only if
$(s\rho,t\rho)\in\delta^{a\rho}$.

If now one has $a=cb$, $b,c\in S^1$, then
$\rho_{b,c}\subseteq\rho_a$. The following statement shows when
the opposite inclusion holds in the case when $S/\rho$ is inverse.

\begin{proposition}
Let $S$ be a semigroup, $b,c\in S^1$. Let also $\rho$ be a
congruence on $S$ such that $S/\rho$ is inverse. Then
$\rho_{b,c}=\rho_{cb}$ if and only if
\begin{equation} \label{main}
uxv=v^{-1}vu\cdot x\cdot vuu^{-1}~ \mbox{for all}~x\in S/\rho,
\end{equation}
where $u=b\rho$ and $v=c\rho$.
\end{proposition}

\begin{proof}
The condition $\rho_{b,c}=\rho_{cb}$ is equivalent to
$\rho_{cb}\subseteq\rho_{b,c}$, which is, in turn, equivalent to
\begin{equation} \label{eq:2}
v\cdot uxv\cdot u=v\cdot uyv\cdot u\Rightarrow uxv=uyv~\mbox{for
all}~x,y\in S/\rho.
\end{equation}
Let us prove that the condition \eqref{eq:2} is equivalent to
\eqref{main}.

Indeed, assume that \eqref{eq:2} holds. Since
\begin{equation*}
vu\cdot x\cdot vu=vu\cdot u^{-1}v^{-1}vu\cdot x\cdot
vuu^{-1}v^{-1} \cdot vu~\mbox{for all}~x\in S/\rho,
\end{equation*}
we have $u\cdot x\cdot v=u\cdot u^{-1}v^{-1}vu\cdot x\cdot
vuu^{-1}v^{-1} \cdot v=v^{-1}vu\cdot x\cdot vuu^{-1}$, due to the
fact that all the idempotents of an inverse semigroup pairwise
commute. That is, the condition \eqref{main} holds.

Assume now that \eqref{main} holds. Suppose that $v\cdot uxv\cdot
u=v\cdot uyv\cdot u$ for some $x,y\in S/\rho$. Then
\begin{equation*}
uxv=v^{-1}vu\cdot x\cdot vuu^{-1}= v^{-1}vu\cdot y\cdot
vuu^{-1}=uyv,
\end{equation*}
whence we obtain \eqref{eq:2}.
\end{proof}

\begin{proposition}
Let $S$ be a semigroup, $a\in S^1$, $\rho\in\mathrm{Cong}(S)$.
Then
\begin{equation*}
S_{a}/\rho_{a}\cong aSa/\rho\cap(aSa\times aSa).
\end{equation*}

\end{proposition}

\begin{proof}
Set a map $\psi:S_{a}\rightarrow aSa/\rho\cap(aSa\times aSa)$ as
follows:
\begin{equation*}
\psi(x)=(axa)\rho~\mbox{for all}~x\in S.
\end{equation*}
It remains to note that $\psi$ is an onto homomorphism and
$\psi\circ\psi^{-1}=\rho_a$.
\end{proof}

\noindent
G.K.: Algebra, Department of Mathematics and Mechanics, Kyiv Taras
Shevchenko University, 64 Volodymyrska st., 01033 Kyiv, UKRAINE,\\
e-mail: {\tt akudr\symbol{64}univ.kiev.ua}
\vspace{0.3cm}

\noindent
V.M.: Algebra, Department of Mathematics and Mechanics, Kyiv Taras
Shevchenko University, 64 Volodymyrska st., 01033 Kyiv, UKRAINE,\\
e-mail: {\tt vmaltcev\symbol{64}univ.kiev.ua}

\end{document}